\documentclass[10pt]{amsart}
     \makeatletter
     \def\section{\@startsection{section}{1}%
     \z@{.7\linespacing\@plus\linespacing}{.5\linespacing}%
     {\bfseries%\normalfont\scshape
     \centering
     }}
     \def\@secnumfont{\bfseries}
     \makeatother
\newtheorem{theorem}{Theorem}[section]
\newtheorem{lemma}[theorem]{Lemma}

\theoremstyle{definition}
\newtheorem{definition}[theorem]{Definition}

\theoremstyle{remark}

\numberwithin{equation}{section}

\newenvironment{namedproof}[1]{\emph{Proof of #1.}}{\hfill$\Box$}

\usepackage{amsmath,amsthm,amssymb,amsbsy}

\usepackage[all]{xy}
\usepackage{chngcntr}
\usepackage{url}
\usepackage{graphicx}

\usepackage{tikz}
\usetikzlibrary{arrows,positioning,calc,shapes.multipart}

\tikzset{
    >=stealth',
    markovnode/.style ={rectangle, rounded corners, draw=black, very thick, text width=6.5em, minimum height=3em, text centered},
    markovAuxNode/.style ={rectangle, rounded corners, draw=white, very thick, text width=6.5em, minimum height=3em, text centered},
    markovnodeSmall/.style ={rectangle, rounded corners, draw=black, very thick, text width=3.5em, minimum height=2em, text centered},
    markovAuxNodeSmall/.style ={rectangle, rounded corners, draw=white, very thick, text width=3.5em, minimum height=2em, text centered},
    markovarrow/.style={->, thick, shorten <=2pt, shorten >=2pt},
    markovarrowdashed/.style={->, dashed, shorten <=2pt, shorten >=2pt}
    }

\counterwithin{figure}{section}

\newcommand{\NN}{\mathbb{N}}

\newcommand{\RR}{\mathbb{R}}

\newcommand{\BBB}{\mathcal{B}}

\newcommand{\EEE}{\mathcal{E}}
\newcommand{\FFF}{\mathcal{F}}
\newcommand{\GGG}{\mathcal{G}}

\newcommand{\PPP}{\mathcal{P}}

\newcommand{\dv}{\,\mathrm{d}}                         % Integrator
\newcommand{\eps}{\varepsilon}                              % Epsilon
\begin{document}

\setlength{\parindent}{0cm}
\setlength{\parskip}{0.5cm}

\title[Revisiting the forward equations]{Revisiting the forward equations for inhomogeneous semi-Markov processes}

\date{\today}

\author[A. Sokol]{Alexander Sokol}

\address{Alexander Sokol: Institute of Mathematics, University of
  Copenhagen, 2100 Copenhagen, Denmark, alexander@math.ku.dk}

\subjclass[2010] {Primary 60K15; Secondary 60J25, 60G55, 35R09}

\keywords{Semi-Markov process, Strong Markov process, Intensity, PDE}

\begin{abstract}
In this paper, we consider a class of inhomogeneous semi-Markov processes directly
based on intensity processes for marked point processes. We show that
this class satisfies the semi-Markov properties defined elsewhere
in the literature. We use the marked point process setting to derive strong upper
bounds on various probabilities for semi-Markov processes. Using these
bounds, we rigorously prove for the case of countably infinite state space that the transition
intensities are right-derivatives of the transition probabilities, and we prove
for the case of finite state space that the transition probabilities satisfy the forward equations,
requiring only right-continuity of the transition intensities in
the time and duration arguments and a boundedness condition. We also show relationships
between several classes of semi-Markov processes considered in the literature, and
we prove an integral representation for the left derivatives of the transition
probabilities in the duration parameter.
\end{abstract}

\maketitle

\noindent

%Mulige udvidelser: Tællelige tilstandsrum.

%Husk: strong Markov giver "| \FFF_T" = "| \sigma(X_T,T)", ikke
%blot \sigma(X_T). Tjek for anvendelser.

%Helwich antager: (Z_m,T_m) homogen Markov.

%Husk : inhomogen implicit, homogen eksplicit

\section{Introduction}

Since the introduction of semi-Markov processes in \cite{Levy1954} and \cite{SmithReg1955}, this
class of stochastic processes have been thoroughly developed and applied in many fields of study. Initially,
the semi-Markov processes studied were homogeneous semi-Markov
processes, see e.g. \cite{Pyke1961,PykeSchaufele1964,FellerSemiMarkov1964,Cinlar1973Survey}.
Intuitively, these are similar to homogeneous Markov processes, except that the intensities of jumps
depend on the amount of time spent by the process in its current state. This duration
dependency is at the center of semi-Markov theory.

Inhomogeneous semi-Markov processes, where the intensity for jumps depend on both time and
the time spent by the process in its current state, have been studied as well,
see e.g. \cite{Hoem1972,JanssenDominicis1984,VassiliouPapadopoulou1992}. Here, the theory
becomes more complex, but also allows for more flexible modeling.

The strength of semi-Markov theory is its ability to model systems where the change of the
state depends on the amount spent in its current state. Many real-world systems
exhibit this type of behaviour. One obvious example is in the context of disability insurance,
where the intensity for recovery depends on the amount of time spent in a disabled state:
Intuitively, if one does not recover within a few years, it is considerably less likely
that one will ever recover. Applications of semi-Markov theory to disability insurance,
and life insurance in general, can be found in e.g.
\cite{ChristiansenMultistate2012,Hoem1972,Helwich2007,BuchardMoellerSchmidt2013}. Semi-Markov models
find applications elsewhere as well, in fields as diverse as for example wind speed modeling,
high frequency finance, tourism movements and credit risk, see \cite{DamicoEtAl2013,DamicoEtAl2013Finance,JianhongEtAl2011,LucasEtAl2006}.

The theoretical literature on semi-Markov processes is considerable, covering
e.g. \cite{Pyke1961,FellerSemiMarkov1964,Hunter1968,Hoem1972,AthreyaEtAl1978,GreenwoodEtAl2004,Helwich2007} and
much more, see the books \cite{JanssenManca2005,JanssenManca2006} for thorough reference lists.
The majority of the literature, however, focuses on the homogeneous case, neglecting the inhomogeneous
case. Furthermore, sometimes the interest in applications take precedence to rigor
and mathematical precision, with e.g. \cite{Hoem1972} stating: \emph{We shall give
a few proofs, but not where the results can be argued by "direct reasoning" (i. e.,
by intuition) \ldots All this means that we shall sweep some interesting mathematical
problems under the rug}.

In this paper, we revisit some fundamental results for semi-Markov processes, namely
the definition of semi-Markov processes, related Markov properties, and the characterization
of intensities and transition probabilities. In particular, we give precise conditions
and a rigorous proof for the transition probabilities to satisfy the forward partial
integro-differential equations, utilizing a direct proof based on estimates on the probabilities of multiple
jumps in a short time interval. The forward equations is one example of a result
related to the citation by \cite{Hoem1972} above: It is intuitively quite clear that the
equations should hold under "sufficient regularity conditions", but obtaining a rigorous
proof of this is challenging.

All of our results cover the inhomogeneous case. Furthermore, we also consider the case of
countably infinite state spaces when possible. Our results not only infuse known results with new and improved proof
methodologies, but also generalize known results. Specifically, our contributions are as follows:

\begin{enumerate}
\item[(1).] We prove a theorem on the relative sizes of various classes of semi-Markov
processes considered in the literature.
\item[(2).] We prove an upper bound for the conditional probability of a semi-Markov
process making two jumps in a small interval, Lemma \ref{lemma:SMGTwoJumpsBound}, and
demonstrate in our proofs that this lemma can be used as a key tool for rigorous proofs
of various analytical properties of semi-Markov processes.
\item[(3).] In the case of a countable state space, we give a precise sufficient condition
on the intensities of a semi-Markov process for being the derivatives of the transition probabilities.
\item[(4).] In the case of a finite state space, we give a precise sufficient condition
on the intensities of a semi-Markov process for satisfying the forward partial
integro-differential equations.
\end{enumerate}

The remainder of the paper is organized as follows. In Section \ref{sec:Definitions}, we review
the various notions of semi-Markov processes considered in the literature, we introduce
a class of semi-Markov processes with extra regularity, and we prove a theorem on the
relative sizes of various classes of semi-Markov processes. In Section \ref{sec:Intensities},
we prove a bound on the probability of a semi-Markov process making two jumps in a short
interval of time, and we give sufficient conditions for the transition intensities to be
the right derivatives of the transition probabilities. In Section \ref{sec:Forward}, we
give sufficient conditions for the transition intensities to satisfy the forward equations.
Finally, in Section \ref{sec:Discussion}, we discuss our results and consider opportunities
for further research. Proofs can be found in Appendix \ref{sec:proofs}.

\section{Weak and strong classes of semi-Markov processes}

\label{sec:Definitions}

Several definitions of the notion of a semi-Markov process exist in the literature,
varying in their strictness. In this section, we review the definitions used in the literature,
and we introduce a class of semi-Markov processes with extra regularity, using
intensity processes for marked point processes. Furthermore, we prove a result on the relative sizes of
various types of semi-Markov process.

We work in the context of a sequence of random variables $(Y_n,S_n)_{n\ge0}$, where $S_0=0$ and $S_n>0$ for
all $n\ge1$ and $Y_n$ takes its values in $E$ for all $n\ge0$, where $(E,\EEE)$ is some measurable space.
We assume that $Y_{n+1}\neq Y_n$ for all $n\ge0$. We then also define $T_n=\sum_{k=0}^n S_k$ for $n\ge0$. In the context of a marked point process,
see \cite{MR2189574}, we think of $(T_n)_{n\ge1}$ as the event times and of $Y_n$ as the marks, with $T_0$ being time zero
and $Y_0$ being the initial mark. We then also define $Z_t = Y_n$ for $T_n\le t < T_{n+1}$, and
define $U_t = t - \sup\{0\le s\le t|Z_s\neq Z_t\}$. We refer to $U$ as the duration process
of $Z$. Intuitively, $U_t$ measures the amount of time spent by $Z_t$ in its current state. Here,
it is implicit that if $T_{n+1}=\infty$ for some $n$, then $Z_t = Y_n$ for all $t\ge T_n$. Note
that using this definition, it is immediate that $Z$ has cadlag sample paths.

One formulation of the semi-Markov property often seen is given directly through the
dependence structure of $(Y_n,S_n)_{n\ge0}$. In \cite{Pyke1961}, the process $Z$ is defined to be a semi-Markov process when
\begin{align}
\label{eq:PykeSemiMarkovChar}
  P(Y_n = j,S_n\le t| Y_0,S_0,Y_1,S_1,\ldots,S_{n-1},Y_{n-1}) = P_{Y_{n-1},j}(t)
\end{align}
for some family $P(i,\cdot)$ of distributions on $E$, where $P_{ij}(t) = P(i,\{j\}\times[0,t])$. Here, then, $P_{ij}(t)$ can
be interpreted as the probability of transitioning from $i$ to $j$, having stayed in the previous state
for a duration of less than or equal to $t$. This definition corresponds to a
time-homogeneous case, and essentially states that $(Y_n,S_n)$ is conditionally
independent of $S_0,Y_0,\ldots,S_{n-1}$ given $Y_{n-1}$. The same definition is also used in \cite{Cinlar1973Survey},
where the family $(P(i,\cdot))_{i\in E}$ is referred to a semi-Markovian kernel. The
papers \cite{DamicoEtAl2013Finance,PykeSchaufele1964} also use variants of this definition.
Furthermore, \cite{Helwich2007} defines $Z$ to be a semi-Markov process if $(Y_n,T_n)_{n\ge0}$
is a Markov process. The majority of authors take the case of finite state space $E$ as their main interest. Some
exceptions to this are \cite{FellerSemiMarkov1964,SmithReg1955,Yackel1966}, who consider
the case of a countable state space, and \cite{VegaAmaya2002}, who allows a general, possibly
uncountable state space.

The other main type of definition considers the process $(Z,U)$. In particular, in \cite{Yackel1966},
$Z$ is said to be a semi-Markov process when $(Z,U)$ is a strong homogeneous Markov process.
This type of definition is also applied in \cite{ChristiansenMultistate2012,BuchardMoellerSchmidt2013},
allowing for time inhomogeneity. In \cite{Hoem1972} the author endeavors to retain the
discussion on an intuitive level, but generally also argues for a definition of semi-Markov processes of the same type.

The definition based on having $(Z,U)$ be a Markov process has several
qualities: It encapsulates the central idea of a semi-Markov process, namely dependence
on only the current state and the duration spent in that state, it is concise and
it is equally coherent for both finite, countable and uncountable state spaces. As stated
in the following definition, we will consider this to be the main defining property
of a semi-Markov process. For sake of tractability, we require that $(Z,U)$ in fact
be a strong Markov process.

\begin{definition}
\label{def:semiMarkovByZU}
Let $(E,\EEE)$ be a measurable space, and let $Z$ be a cadlag process taking its values
in $E$. We say that $Z$ is an \emph{inhomogeneous semi-Markov process} if $(Z,U)$ is
an inhomogeneous strong Markov process on $E\times\RR_+$.
\end{definition}

By itself, however, processes such as those given in Definition \ref{def:semiMarkovByZU}
are not sufficiently regular to allow for a rich mathematical theory. We next introduce a
more regular type of semi-Markov processes, based on intensity processes
for marked point processes. Our claim is that this class of processes is
amenable to discuss various types of regularity and will provide a sound framework for
rigorous development of theory.

In the following, we assume given for each $i,j\in E$ with $i\neq j$ a
measurable mapping $q_{ij}:\RR_+^2\to\RR_+$. We also define $q_i:\RR_+^2\to\RR_+$ for $i\in E$
by
\begin{align}
	q_i(t,u)=\sum_{j\neq i}q_{ij}(t,u),
\end{align}
and we assume throughout that
\begin{align}
\label{eq:FiniteIntensityBound}
	\sup_{(t,u)\in K} \sup_{i\in E} q_i(t,u)<\infty
\end{align}
for all compact subsets $K$ of $\RR_+^2$. This
will in particular ensure the absence of explosion for all processes under
consideration. The requirement (\ref{eq:FiniteIntensityBound}) may appear
strict at first sight, but really states little more than that all intensities for
making jumps are bounded simultaneously on finite intervals of time. In the case
where $E$ is finite, the supremum over $E$ is of course of no impact. We furthermore define
$q_{ii}(t,u)=-q_i(t,u)$, and we let $Q(t,u)$ denote the $E\times E$ matrix
whose entry for the $i$'th row and $j$'th column is $q_{ij}(t,u)$.

\begin{definition}
\label{def:SemiMarkov}
We say that a marked point process process $Z$ with countable state space $E$ is an \emph{inhomogeneous semi-Markov
process with intensities} if there exists a family of mappings $q_{ij}:\RR_+^2\to\RR_+$ for $i\neq j$
satisfying (\ref{eq:FiniteIntensityBound}) such that $Z$ has intensity process
$\lambda$ given by $\lambda_t(i) = 1_{(Z_{t-}\neq i)}q_{Z_{t-}i}(t,U_{t-})$
for $i\in E$, where we assume that each $q_i$ is Lebesgue integrable on compact subsets
of $\RR_+^2$.
\end{definition}

Note that by Corollary 4.4.4 of \cite{MR2189574}, the bound (\ref{eq:FiniteIntensityBound})
ensures that explosion does not occur and thus Definition \ref{def:SemiMarkov} is not vacuous.
In Theorem \ref{theorem:SemiMarkovYieldsStrongMP}, we show that the processes
defined in Definition \ref{def:SemiMarkov} are in fact semi-Markov processes
in the sense of Definition \ref{def:semiMarkovByZU}. The requirement that $q_i$ is Lebesgue integrable
on compact subsets of $\RR_+^2$ is made to ensure that the resulting intensity
process $(t,i)\mapsto \lambda_t(i)$ in fact is integrable and thus allows for
a corresponding compensator process. The interpretation of Definition \ref{def:SemiMarkov} is
straightforward: When the process $Z$ has been in state $i$
for a duration $u$ at time $t$, the intensity for making a jump to a state
$j\neq i$ is $q_{ij}(t,u)$.

We now show a result, Theorem \ref{theorem:SemiMarkovVariantRelations}, about the relative
strengths of various notions of semi-Markov processes. The result that 
Definition \ref{def:SemiMarkov} implies Definition \ref{def:semiMarkovByZU} is
sufficiently important to us to be stated as a separate theorem.

\begin{theorem}
\label{theorem:SemiMarkovYieldsStrongMP}
Assume that $Z$ is an inhomogeneous semi-Markov process with intensities.
Then $(Z,U)$ is an inhomogeneous strong Markov process on $E\times\RR_+$.
\end{theorem}

\begin{theorem}
\label{theorem:SemiMarkovVariantRelations}
Let $Z$ be a cadlag stochastic process with countable state space $E$, let $(T_n)$
be the jump times of $Z$, and let $U$ be the duration process of $Z$. Consider the following
four statements:
\begin{enumerate}
\item[(1).] The process $Z$ is an inhomogeneous semi-Markov process with intensities.
\item[(2).] The process $Z$ is an inhomogeneous semi-Markov process.
\item[(3).] $(Y_n,S_n)$ is conditionally independent of $S_0,Y_0,\ldots,S_{n-1}$ given $Y_{n-1}$.
\item[(4).] $(Y_n,T_n)_{n\ge0}$ is a discrete-time Markov chain.
\end{enumerate}
Here, (1) refers to Definition \ref{def:SemiMarkov}, and (2)
refers to definition Definition \ref{def:semiMarkovByZU}. It then holds that $(1)$ implies $(2)$, $(2)$ implies $(3)$,
and $(3)$ implies $(4)$.
\end{theorem}

In the next section, we consider results characterizing the distribution of
a semi-Markov process with intensities, in particular characterizing the transition probabilities of
the process.

\section{Intensities as derivatives of transition probabilities}

\label{sec:Intensities}

Our objective in this section will be to rigorously obtain sufficient regularity
criteria for the intensity functions $q_{ij}$ to ensure that they can be
obtained as the right derivatives of the transition probabilities of the
semi-Markov process. Furthermore, we will show that under the same
regularity conditions, the convergence of the difference quotients occur under
the presence of a dominating bound. This will be essential for our later
proof that the transition probabilities satisfy the forward equations.

In the following, let $Z$ be a semi-Markov process with intensities $q_{ij}$ for
$i\neq j$. As in the previous section, we assume that the intensities
satisfy (\ref{eq:FiniteIntensityBound}). As the process $(Z,U)$ is a Markov process
by Theorem \ref{theorem:SemiMarkovYieldsStrongMP}, we may associate to it
a family of transition probabilities $P_{s,t}((i,u),\cdot)$ on $E\times \RR_+$ such that
\begin{align}
  P((Z_t,U_t)\in A\times B \mid \FFF_s) &= P_{s,t}(Z_s,U_s,A\times B),
\end{align}
where $A\subseteq E$ and $B\in\BBB_+$, with $\BBB_+$ denoting the Borel $\sigma$-algebra
on $\RR_+$. Consistently with the notation outlined in \cite{Hoem1972} and also used in
\cite{BuchardMoellerSchmidt2013}, we define, for $i,j\in E$ and $s,t,u,d\ge0$
with $s\le t$,
\begin{align}
p_{ij}(s,t,u,A)&=P_{s,t}(i,u,\{j\}\times A), \\
p_{ij}(s,t,u,d)&=p_{ij}(s,t,u,[0,d]),\\
p_{ij}(s,t,u)&= p_{ij}(s,t,u,\RR_+),
\end{align}
where $A\in\BBB_+$. We then think of $p_{ij}(s,t,u,A)$ as the probability of transitioning from $i$ to $j$
from time $s$ to time $t$ when the duration at time $s$ is $u$, and requiring that the duration
at time $t$ is in $A$. We then also similarly think of $p_{ij}(s,t,u)$ as the probability
of transitioning from $i$ to $j$ from time $s$ to time $t$ when the duration at time $s$ is $u$.
Finally, as a function of $d$, $p_{ij}(s,t,u,d)$ is the cumulative mass function for the
measure $A\mapsto p_{ij}(s,t,u,A)$.

Our objective is to show that the intensities $q_{ij}$ are right derivatives of transition
probabilities. To this end, we first state two lemmas of independent interest to be used in
the proof. For use in the following, we define
\begin{align}
N_t &= \sum_{0<s\le t}1_{(Z_{s-}\neq Z_s)},\label{eq:SMGDefAllJumps}
\end{align}
the counting process for counting the jumps of $Z$. Also, when $Q$ is some
$E\times E$ matrix, we let $\|Q\|_\infty$ denote the supremum norm of the matrix. Note
that even when $E$ is infinite, it always holds for $t,u\ge0$ that
\begin{align}
\label{eq:FiniteMatrixIntensityBound}
  \sup_{(t,u)\in K}\|Q(t,u)\|_\infty
  &=\sup_{(t,u)\in K}\sup_{i,j\in E}|q_{ij}(t,u)|
   =\sup_{(t,u)\in K}\sup_{i\in E}q_i(t,u)<\infty,
\end{align}
for compact subsets $K$ of $\RR_+^2$, due to the condition (\ref{eq:FiniteIntensityBound}).
Finally, we say that the intensities $(q_{ij})$ are
right-continuous if it holds for all $i\neq j$ that $q_{ij}(t+h,u+k)$ tends to $q_{ij}(t,u)$
whenever $h$ and $k$ tends to zero from above.

\begin{lemma}
\label{lemma:SMGTwoJumpsBound}
It holds that 
\begin{align}
  &P(N_{t+h}-N_t\ge 2|Z_t = i,U_t=u)\notag\\
  &\le \int_t^{t+h} \int_s^{t+h} \|Q(s,u+s-t)\|_\infty\|Q(v,v-s)\|_\infty \dv v\dv s.
\end{align}
\end{lemma}

\begin{lemma}
\label{lemma:SMGNoTwoJumps}
It holds that
\begin{align}
\lim_{h\to0}\frac{1}{h} P(N_{t+h}-N_t\ge 2|Z_t = i,U_t=u)&=0.
\end{align}
\end{lemma}

Lemma \ref{lemma:SMGNoTwoJumps} essentially shows that for semi-Markov processes
satisfying the bound (\ref{eq:FiniteIntensityBound}), the probability of having two jumps in a small time
interval tends to zero at a faster than linear rate. This is an important regularity
property which is one of the key properties allowing us to prove our results on semi-Markov
processes. The lemma follows immediately from Lemma \ref{lemma:SMGTwoJumpsBound}.
Lemma \ref{lemma:SMGTwoJumpsBound} is a stronger bound, allowing for a variety
of derived bounds on probabilities for semi-Markov processes.

We are now ready to state our two main results of this section. Theorem \ref{theorem:IntensityIsTransProbLimit}
shows that the intensities are limits of transition probabilities, in accordance
with our intuitive understanding of intensities. Theorem \ref{theorem:IntensityIsTransProbLimitBounded}
shows that the convergence in Theorem \ref{theorem:IntensityIsTransProbLimit} occurs
under the presence of a dominating bound. This latter result is used in
the proof of Theorem \ref{theorem:IntensityIsTransProbLimit} and will also
be essential for our rigorous proof of the forward equations in the next section. In
the statement of the theorems, we let $I$ denote the $E\times E$ identity matrix.

\begin{theorem}
\label{theorem:IntensityIsTransProbLimitBounded}
If the intensities are right-continuous, there exists a measurable mapping $C:\RR_+^2\to\RR_+$, bounded on compacts,
such that there is $\eps>0$ with the property that for any $t,u\ge0$ and $0\le h\le \eps$, we have
\begin{align}
\label{eq:IntensityTransProbLimitBound}
  \sup_{i\in E}\sum_{j\in E}\left|\frac{1}{h}(p_{ij}(t,t+h,u)-I_{ij})\right|
  \le C(t,u).
\end{align}
\end{theorem}

\begin{theorem}
\label{theorem:IntensityIsTransProbLimit}
If the intensities are right-continuous, it holds that
\begin{align}
\label{eq:IntensityTransProbLimit}
  \lim_{h\to0} \frac{1}{h}(P(t,t+h,u)-I) = Q(t,u).
\end{align}
\end{theorem}

\section{The forward equations}

\label{sec:Forward}

Next, we turn to the characterization of transition probabilities. For homogeneous
Markov processes, a closed-form expression for the transition probabilities is available to us through
the matrix exponential, see e.g. Section 2.8 of \cite{MR1600720}, and for inhomogeneous Markov processes,
the transition probabilities can be characterized through two different systems
of multidimensional ODEs, namely the forward and backward equations, as proven
in \cite{Feller1940}. In the semi-Markov case, the transition probabilities
also satisfy forward and backward systems of equations. In the case of the forward
equations, these take the form of a system of partial integro-differential equations.
The existence of these equations for semi-Markov processes is
well known. That the transition probabilities satisfy both the forward and backward
equations is proved very rigorously in \cite{Helwich2007} for the case of a homogeneous semi-Markov process. A variant of
the forward equations for the inhomogeneous case is stated in \cite{BuchardMoellerSchmidt2013},
where it forms the basis for numerical calculations of cashflows in life insurance.

In this section, we use the results from Section \ref{sec:Intensities} to give
a direct and rigorous proof of a sufficient condition on the intensities for
the transition probabilities to satisfy the forward equations. We will also
show that this set of equations can be rewritten as a set of ordinary
integro-differential equations in two different ways.

Throughout this section, we assume that \emph{$E$ is finite and that the transition
probabilities are right-continuous}. The condition (\ref{eq:FiniteIntensityBound})
also remains in force, but is now equivalent to the simpler condition
\begin{align}
\label{eq:FiniteIntensityBoundFiniteE}
	\sup_{(t,u)\in K} q_i(t,u)<\infty.
\end{align}
%Bemærk: Cadlag er begrænsede på kompakte intervaller, men det samme gælder ikke
%hvis kun højrekontinuitet (og ikke grænser fra venstre) kræves.
Before proving that the transition probabilities satisfy the forward equations,
we require two lemmas of independent interest. The first lemma states that in an asymptotic sense, no double jumps
back and forth are made over a small time period by a semi-Markov chain. Here,
we define $p_{ij}(s,t,u,d-) = \lim_{h\to0^+}p_{ij}(s,t,u,d-h)$ for $d>0$.
	
\begin{lemma}
\label{lemma:SMGNoQuickCycles}
It holds that 
\begin{align}
\label{eq:NoQuickCyclesLimit}
\lim_{h\to0}\frac{1}{h} p_{ii}(t,t+h,u,(u+h)-)=0.
\end{align}
Furthermore, there exists a measurable mapping $C:\RR_+^2\to\RR_+$, bounded on
compacts, such that there is $\eps>0$ with the property that for any $t,u\ge0$
and $0\le h\le \eps$, we have
\begin{align}
\label{eq:NoQuickCyclesBound}
  \frac{1}{h} p_{ii}(t,t+h,u,(u+h)-) \le C(t,u).
\end{align}
\end{lemma}

Note that the result in Lemma \ref{lemma:SMGNoQuickCycles} cannot be extended
to cover $p_{ii}(t,t+h,u,u+h)$ instead of $p_{ii}(t,t+h,u,(u+h)-)$, since
$p_{ii}(t,t+h,u,\cdot)$ is concentrated on $[0,u+h]$ and so it holds that
\begin{align}
\lim_{h\to0}\frac{1}{h} p_{ii}(t,t+h,u,u+h)
&=\lim_{h\to0}\frac{1}{h} p_{ii}(t,t+h,u)\notag\\
&=\lim_{h\to0}\frac{1}{h} (p_{ii}(t,t+h,u)-1)+\frac{1}{h},\notag
\end{align}
which is infinite, according to Theorem \ref{theorem:IntensityIsTransProbLimit}. Essentially, this shows that in an asymptotic sense,
the only mass for $p_{ii}(t,t+h,u,\cdot)$ is a point mass in $u+h$, corresponding
to no jumps being made when transitioning from $i$ to $i$ from time $t$
to time $t+h$.

We are now ready to begin work on proving the differential properties of
the transition probabilities. We begin by showing that the transition probabilities
are differentiable from the left in the final parameter, and we derive an explicit expression for the
derivative. With this lemma in hand, we will be able to prove our main result
on the forward equations.

\begin{lemma}
\label{lemma:pIsDurationDiff}
Let $i,j\in E$, let $0\le s\le t$ and let $u\ge0$. It holds that
the mapping $d\mapsto p_{ij}(s,t,u,d)$ is differentiable from the left for $d>0$. Furthermore,
for $i\neq j$ and $d>t-s$, the derivative is zero, for $i\neq j$ and $d\le s-t$,
the derivative is 
\begin{align}
  \frac{\dv p_{ij}}{\dv d}(s,t,u,d)
  &=\sum_{k\neq j} \int_0^{u+t-d-s} a_{kj}(t,d,u,v) p_{ik}(s,t-d,u,\dv v).\label{eq:TransProbDurationDiffNonDiag}
\end{align}
where
\begin{align}
  a_{kj}(t,d,u,v) &= q_{kj}(t-d,v)\exp\left(-\int_{t-d}^t q_j(r,r-(t-d))\dv r\right).\label{eq:TransProbDurationDiffNonDiagHelper}
\end{align}
\end{lemma}

The intuitive explanation for the formula (\ref{eq:TransProbDurationDiffNonDiag}) is as
follows. The derivative is the limit as $h$ tends to zero from above of
\begin{align}
  \frac{1}{h}\left(p_{ij}(s,t,u,d) - p_{ij}(s,t,u,d-h)\right)
  &=\frac{1}{h}p_{ij}(s,t,u,(d-h,d])\notag\\
  &=\frac{1}{h}P(Z_t = j,d-h<U_t\le d | Z_s = i,U_s = u).\label{eq:durationDiffApprox}
\end{align}
Now, having $d-h<U_t\le d$ is equivalent to having a jump being made in the
time interval $[t-d,t-d+h)$ and having no jumps made in the time
interval $[t-d+h,t]$. Conditioning on $Z_{t-d} = k$ for $k\in E$ yields that the the differential
quotient (\ref{eq:durationDiffApprox}) approximately is the sum over the states $k$,
and for each state, we sum the density of transitioning from $k$ to $j$ immediately
after time $t-d$ and then remaining in state $j$ from time $t-d$ to time $t$
with duration zero at time $t-d$. Furthermore, this is weighted with the probability
of transitioning from state $i$ to state $k$ from time $s$ to time $t-d$ with duration $u$ at time $s$.
The term corresponding to $k=j$ vanishes, as being in state $j$ at time $t-d$
would indicate an ultimate duration at time $t$ greater than $d$. In
(\ref{eq:TransProbDurationDiffNonDiag}), the integration with respect to $p_{ik}(s,t-d,u,\dv v)$
represents conditioning on the state at time $t-d$, the left factor in
(\ref{eq:TransProbDurationDiffNonDiagHelper}) corresponds to the conditional density of making a
jump from $k$ to $j$ immediately after time $t-d$, and the right factor in
(\ref{eq:TransProbDurationDiffNonDiagHelper}) corresponds to the probability of remaining
in state $j$ from time $t-d$ to time $t$ with duration zero at time $t-d$.

We are now ready to prove that the transition probabilities satisfy the
forward partial integro-differential equations.

\begin{theorem}
\label{theorem:SMPForward}
Fix $s\ge0$ and $u\ge0$. It holds for $t\ge s$ and $d> 0$ that
$p_{ij}(s,t,u,d)$ is differentiable in $t$ from the right, and the derivative
is
\begin{align}
\label{eq:SMGForward}
\frac{\partial p_{ij}}{\partial t}(s,t,u,d)
 &=\sum_{k\neq j} \int_0^{u+t-s} q_{kj}(t,v) p_{ik}(s,t,u,\dv v)\notag\\
 &+\int_0^d q_{jj}(t,v) p_{ij}(s,t,u,\dv v)-\frac{\partial p_{ij}}{\partial d}(s,t,u,d),
\end{align}
where the partial derivative with respect to $d$ is the derivative from the left.
\end{theorem}

Theorem \ref{theorem:SMPForward} yields the semi-Markovian analogue of the
Kolmogorov forward equations for Markov processes. As an immediate corollary of
Theorem \ref{theorem:SMPForward}, we may also derive a system of ordinary
integro-differential equations for the transition probabilities.
Applying the chain rule, it holds that
\begin{align}
\frac{\partial}{\partial t}p_{ij}(s,t,u,d+t-s)
&=\frac{\partial p_{ij}}{\partial t}(s,t,u,d+t-s)
 +\frac{\partial p_{ij}}{\partial d}(s,t,u,d+t-s),\label{eq:ForwardODE1}
\end{align}
and applying Theorem \ref{theorem:SMPForward} in (\ref{eq:ForwardODE1}), we obtain
\begin{align}
\label{eq:SMGForwardODE}
\frac{\partial}{\partial t}p_{ij}(s,t,u,d+t-s)
&=\sum_{k\neq j} \int_0^{u+t-s} q_{kj}(t,v) p_{ik}(s,t,u,\dv v)\notag\\
&+\int_0^{d+t-s} q_{jj}(t,v) p_{ij}(s,t,u,\dv v).
\end{align}
Note that in (\ref{eq:ForwardODE1}) and (\ref{eq:SMGForwardODE}), the term
\begin{align}
\frac{\partial}{\partial t}p_{ij}(s,t,u,d+t-s)
\end{align}
refers to differentiation of the composite mapping $p_{ij}(s,t,u,d+t-s)$, and not to
the derivative of $p_{ij}(s,t,u,d)$ evaluated in $(s,t,u,d+t-s)$.  Also note that in the
case of a Markov process, where $q_{ij}(t,v)=q_{ij}(t)$ for some
$q_{ij}:\RR_+\to\RR_+$, we may let $u=d$ in (\ref{eq:SMGForwardODE}) and immediately
recover the classical Kolmogorov forward equation from \cite{Feller1940}.

Furthermore, we may also insert our expression (\ref{eq:TransProbDurationDiffNonDiag})
for the derivative in $d$ of $p_{ij}(s,t,u,d)$ directly into (\ref{eq:SMGForward}) to obtain
\begin{align}
\label{eq:SMGForwardODEByDiffd}
\frac{\partial p_{ij}}{\partial t}(s,t,u,d)
 &=\sum_{k\neq j} \int_0^{u+t-s} q_{kj}(t,v) p_{ik}(s,t,u,\dv v)
 +\int_0^d q_{jj}(t,v) p_{ij}(s,t,u,\dv v)\notag\\
 &-\sum_{k\neq j} \int_0^{u+t-d-s} a_{kj}(t,d,u,v) p_{ik}(s,t-d,u,\dv v).
\end{align}
for the case $i\neq j$ and $d\le t-s$, with a similar result holding for the diagonal case.

\section{Discussion}

\label{sec:Discussion}

In this article, we have considered a class of regular semi-Markov processes defined
through intensities for a marked point process. We have proven relationships between
different notions of semi-Markov processes, and we have shown how to obtain sufficient
conditions on the intensities for having the intensities be the right derivatives of
the transition probabilities. Finally, we have given sufficient conditions on the
intensities for having the transition probabilities satisfy the forward equations,
and we have stated the latter equations in three different variants: in one way
as a system of partial integro-differential equations, and in two ways as a system
of ordinary integro-differential equations. In the course of this, we have also proved
a formula for the derivative of the transition probabilities with respect to the
duration $d$.

Our main purpose and focus has been to develop a framework and tools for rigorous
analysis of semi-Markov processes. Our results indicate
that analysis of semi-Markov processes centered around conditioning arguments
for single jumps, combined with applications of bounds for probabilities of
several jumps such as Lemma \ref{lemma:SMGTwoJumpsBound} and Lemma \ref{lemma:SMGNoTwoJumps},
is a fruitful proof methodology.

We feel that there is ample room for improvement of our results and for continued
study of semi-Markov processes. Some subjects
for further study of semi-Markov processes which we find particularly interesting are as follows:
\begin{enumerate}
\item[(1).] Assuming sufficient regularity conditions, the transition probabilities
also satisfy a different system of equations, the backward partial differential
equations, see e.g. \cite{BuchardMoellerSchmidt2013,ChristiansenMultistate2012}, and
see \cite{Helwich2007} for a rigorous proof in the homogeneous case. It is of interest
to develop rigorous arguments for sufficient conditions on the intensities ensuring
the validity of the backward equations.
\item[(2).] In terms of the forward equations, it is of interest to weaken
the regularity conditions required for the intensities, and to consider the case
of countably infinite state spaces. In the latter scenario, several interchanges of
limits and summation would have to be argued for separately.
\item[(3).] The forward and backward equations can be used for numerical evaluation
of various expressions in semi-Markov models, see e.g. \cite{BuchardMoellerSchmidt2013}.
In terms of numerical studies, it would be of interest to compare the computational
efficiency and accuracy of numerical solutions of these equation systems.
\end{enumerate}

We hope that our efforts in this paper will inspire further study of the field
of semi-Markov processes.

\appendix

\section{Proofs}

\label{sec:proofs}

\subsection{Proofs for Section \ref{sec:Definitions}} \hfill\break

\begin{namedproof}{Theorem \ref{theorem:SemiMarkovYieldsStrongMP}}
We may assume without loss of generality that $Z$ has deterministic initial state $y_0$. We
begin by proving that the process $(Z,U)$ is a piecewise deterministic Markov process obtained from a marked point process satisfying
the sufficient criteria of Theorem 7.3.1 of \cite{MR2189574}. In order to do so,
we need to specify both the underlying marked point process as well as the transform
functions $\phi_{s,t}$ appearing in that theorem.

Define $G = E\times \RR_+$. Letting $E$ be endowed with the discrete topology and letting
$\RR_+$ be endowed with its usual topology, we endow $G$ with the product topology. With
$\GGG$ denoting the corresponding Borel $\sigma$-algebra on $G$, it holds that $\GGG$
is countably generated and contains all singletons, as is required by the results
of Section 7.3 of \cite{MR2189574}. Next, let $\tilde{Z}_t = (Z_t,0)$ and define $\tilde{T}_n = T_n$ and $\tilde{Y}_n = (Y_n,0)$.
We then obtain that $\tilde{Z}$ is a marked point process with mark space $(G,\GGG)$.
Furthermore, for $0\le s\le t$ and with $\tilde{y}=(y,u)$, define $\phi_{s,t}:G\to G$ by
$\phi_{s,t}(\tilde{y}) = (y,t-s)$. Also define
$\tilde{T}_{\langle t\rangle}=\sup\{\tilde{T}_n|\tilde{T}_n\le t\}$ and define $\tilde{Y}_{\langle t\rangle} = \tilde{Y}_{\tilde{T}_{\langle t\rangle}}$.
We then obtain
\begin{align}
	(Z_t,U_t) &= (Y_{\langle t\rangle},t-T_{\langle t\rangle})
	           = \phi_{T_{\langle t\rangle},t}(\tilde{Y}_{\langle t\rangle}).
\end{align}
This shows that $(Z,U)$ can be obtained as a piecewise deterministic
process as in Theorem 7.3.1 of \cite{MR2189574}, with underlying marked point
process $\tilde{Z}$. It remains to check the requirements on $\tilde{Z}$ of
that theorem. In the notation of \cite{MR2189574}, the distribution of the marked point
process $\tilde{Z}$ satisfies
\begin{align}
\bar{P}^{(n)}_{\tilde{z}_n}(t) &= \exp\left(-\int_{t_n}^t  q_{y_n}(s,s-t_n)\dv s\right),\\
\pi^{(n)}_{\tilde{z}_n,t}(\{(j,0)\}) &= \frac{q_{y_nj}(t,t-t_n)}{q_{y_n}(t,t-t_n)},\label{eq:PiUnderlyingMPP}
\end{align}
for $t\ge0$, $n\ge1$ and $j\in E$ with $(j,0)\neq \tilde{y}_n$ and with (\ref{eq:PiUnderlyingMPP})
the latter being defined to be zero when the denominator is zero,
the value of $\pi^{(n)}_{\tilde{z}_n,t}(\{(j,0)\})$ in this case is irrelevant for the
distribution of the marked point process. Here, $\tilde{z}_n=(t_1,\tilde{y}_1,\ldots,t_n,\tilde{y}_n)$
and $\tilde{y}_n = (y_n,0)$. Also, similar expressions are obtained
for the distribution of $(\tilde{T}_1,\tilde{Y}_1)$. With $\tilde{y} = (y,u)\in G$, now define
$\tilde{q}_t(\tilde{y}) = q_y(t,u)$. Noting that $\phi_{0,t}(y,u) = (y,t)$, we obtain
\begin{align}
  q_{y_0}(s,s) &= \tilde{q}_s(y_0,s)= \tilde{q}_s(\phi_{0,s}(\tilde{y}_0)), \\
 q_{y_n}(s,s-t_n) &= \tilde{q}_s(y_n,s-t_n) = \tilde{q}_s(\phi_{t_n,s}(\tilde{y}_n)),
\end{align}
so that with $t\ge0$, $(y,u)\in G$ and $C\in\GGG$ and
\begin{align}
	r_t((y,u),C)
	&=\sum_{(j,v)\in C}1_{(j\neq y,v=0)}\frac{q_{yj}(t,u)}{q_{y}(t,u)},
\end{align}
we obtain
\begin{align}
\bar{P}^{(n)}_{\tilde{z}_n}(t) &= \exp\left(-\int_{t_n}^t \tilde{q}_s(\phi_{t_n,s}(\tilde{y}_n))\dv s\right),\\
\pi^{(n)}_{z_n,t}(C) &= r_t(\phi_{t_n,t}(\tilde{y}_n),C).
\end{align}
Finally, we have that
\begin{align}
	\int_t^{t+h} \tilde{q}_s(\phi_{s,t}(\tilde{y}))\dv s
    &=\int_t^{t+h} q_y(s,t-s)\dv s
\end{align}
is finite for small $h$ by assumption, so that (i) of Theorem 7.3.1 of \cite{MR2189574}
is satisfied. As the remaining requirements (ii) and (iii) are trivially satisfied,
we may invoke the theorem and obtain that $(Z,U)$ is a piecewise deterministic Markov
process obtained from a marked point process satisfying the sufficient criteria of
Theorem 7.3.1 of \cite{MR2189574}. Theorem 7.5.1 of \cite{MR2189574}
then furthermore shows that $(Z,U)$ is a strong inhomogeneous Markov process.
\end{namedproof}

\begin{lemma}
\label{lemma:SMGSimplePathProperties}
For all $n\ge1$, it holds that when $T_n<\infty$, we have $Z_{T_n} = Y_n$ and $Z_{T_n-}=Z_{T_{n-1}}$.
Furthermore, for $T_n\le t<T_{n+1}$, it holds that $U_t = t - T_n$. In particular, $U_{T_n-} = S_n$,
with $S_n = T_n - T_{n-1}$.
\end{lemma}

\begin{namedproof}{Lemma \ref{lemma:SMGSimplePathProperties}}
The result on the values of $Z_{T_n}$ and $Z_{T_{n-1}}$ follow immediately from
the pathwise definition of $Z$ in terms of $T_n$ and $Y_n$. As regards the claims about
$U$, let $T_n\le t<T_{n+1}$. It then holds that $Z$ is constant on $[T_n,t]$, while
it holds that $Z_{T_n-} = Y_{n-1} \neq Y_n =Z_{T_n}$, so $Z$ changes its value at $T_n$. Therefore,
we have
\begin{align}
	  \sup\{0\le s\le t | Z_s\neq Z_t\}
	  &=T_n,
\end{align}
and so $U_{T_n} = t- T_n$. In particular, it follows that
\begin{align}
	U_{T_n-}
	&=\lim_{h\to 0^+} U_{T_n-h}
	 =\lim_{h\to 0^+} (T_n-h) - T_{n-1}
	 =T_n-T_{n-1} = S_n,
\end{align}
as required.
\end{namedproof}

\begin{namedproof}{Theorem \ref{theorem:SemiMarkovVariantRelations}}
That $(1)$ implies $(2)$ is the content of Theorem \ref{theorem:SemiMarkovYieldsStrongMP}. To prove
that $(2)$ implies $(3)$, assume that $(Z,U)$ is an inhomogeneous Markov process. Note
that as $S_0=0$, it suffices to prove that $(Y_n,S_n)$ is conditionally independent of $Y_0,S_1\ldots,Y_{n-2},S_{n-1}$
given $Y_{n-1}$. Noting that we always have $U_{T_{n-1}} = 0$, we may apply Lemma \ref{lemma:SMGSimplePathProperties}
and the strong Markov property at the stopping time $T_{n-1}$, obtaining
\begin{align}
  & P(Y_n=j, S_n \le t|Y_0=i_0,S_1=s_1,\ldots,S_{n-1} = s_{n-1},Y_{n-1}=i_{n-1})\notag\\
  &=P(Z_{T_n}=j, U_{T_n-} \le t|Y_0=i_0,S_1=s_1,\ldots,S_{n-1} = s_{n-1},Y_{n-1}=i_{n-1})\notag\\
  &=P(Z_{T_n}=j, U_{T_n-} \le t|U_{T_{n-1}} = 0,Z_{T_{n-1}}=i_{n-1})\notag\\
  &=P(Y_n=j, S_n \le t|Y_{n-1}=i_{n-1}),
\end{align}
which shows the desired conditional independence statement. As regards the proof that $(3)$ implies $(4)$, we note that with
$s_k = t_k-t_{k-1}$ for $1\le k\le n-1$, it holds with $a=t-t_{n-1}$ that
\begin{align}
  & P(Y_n=j, T_n \le t|Y_0=i_0,T_1=t_1,\ldots,T_{n-1} = t_{n-1},Y_{n-1}=i_{n-1})\notag\\
  &=P(Y_n=j, S_n \le a|Y_0=i_0,S_1=s_1,\ldots,S_{n-1} = s_{n-1},Y_{n-1}=i_{n-1})\notag\\
  &=P(Y_n=j, S_n \le t-t_{n-1}|Y_{n-1}=i_{n-1}).
\end{align}
As this does not depend on $t_1,\ldots,t_{n-2}$, and neither on $i_0,\ldots,i_{n-2}$,
we conclude that $(Y_n,T_n)$ is conditionally independent of
$Y_0,T_0,\ldots,Y_{n-2},T_{n-2}$ given $Y_{n-1}$ and $T_{n-1}$, yielding the desired
discrete-time Markov property.
\end{namedproof}

\subsection{Proofs for Section \ref{sec:Intensities}}\hfill\break

Before proving the main results of Section \ref{sec:Intensities}, we first
prove a series of lemmas. To this end, we fix some notation and recall
some facts about hazard functions and hazard measures. First, for $t\ge0$
and any stopping time $S$ with respect to the filtration $(\FFF_t)$ induced
by $Z$, let $T(S)$ denote the first jump strictly after time $S$. Note that as
we have $T(S)=\inf\{T_n|T_n > S\}=\inf (T_n)_{(T_n>S)}$ and $(T_n>S)\in\FFF_{T_n}$, see
e.g. Section I.1b of \cite{MR1943877}, it holds that $T(S)$ is a stopping time for any $S$.
Also, recall that for a distribution $\mu$ on $\RR_+\cup\{\infty\}$ whose restriction
to $\RR_+$ has density $f$ with respect to the Lebesgue
measure, we may define the hazard function $h:\RR_+\to\RR_+$ of $\mu$ by $h(t) = f(t) / \mu((t,\infty])$
whenever the denominator is nonzero, otherwise we let $h(t)=0$, and it then holds that
\begin{align}
  \mu((t,\infty]) &= \exp\left(-\int_0^t h(s)\dv s\right),\\
  f(t) &= h(t)\exp\left(-\int_0^t h(s)\dv s\right).
\end{align}
For distributions which are not absolutely continuous, the more general construct of
a hazard measure can be applied, see Section 4.1 of \cite{MR2189574}. The following
four lemmas yield distributions and hazard functions for various variables
related to the distribution of $Z$. In the following, we let $D(E)$ denote
the set of cadlag mappings from $\RR_+\to E$, and generally apply the notation
of \cite{MR2189574} when considering distributional properties of marked
point processes.

%See (4.47) of Jacobsen for the source of the main formula for the hazard measure.
\begin{lemma}
\label{lemma:ConditionalEventHazardSMP}
The conditional distribution of $T_{n+1}$ given $T_0,Y_0,\ldots,T_n,Y_n$ has hazard function
$h:\RR_+\to\RR_+$ given by $h(t) = 1_{(t>T_n)}q_{Y_n}(t,t-T_n)$.
\end{lemma}
\begin{namedproof}{\ref{lemma:ConditionalEventHazardSMP}}
With $z_n=(t_0,y_0,\ldots,t_n,y_n)$, let $P^{(n)}_{z_n}$ denote the conditional distribution of $T_{n+1}$ given
$(T_0,Y_0,\ldots,T_n,Y_n) = z_n$, and let $\Lambda:D(E)\to D(E)$ denote the total compensator
of $Z$. From the explicit representation of the
compensator of counting processes given in Section 4.3 of \cite{MR2189574},
we obtain that the hazard measure of the conditional distribution is concentrated
on $(t_n,\infty)$ and given by
\begin{align}
\nu^{(n)}_{z_n}((t_n,t]) &= \Lambda_t(z)-\Lambda_{t_n}(z)
=\sum_{j\neq y_n}\int_{t_n}^t q_{y_nj}(s,u_{s-})\dv s\notag\\
 &=\int_{t_n}^t q_{y_n}(s,s-t_n)\dv s,
\end{align}
see (4.47) of \cite{MR2189574}, for paths $z$ such that only $n$ jumps are made on $[0,t)$ and corresponding to
having jump times $t_1,\ldots,t_n$ with destination states $y_1,\ldots,y_n$. As the hazard
measure is absolutely continuous, we obtain that the hazard function for the distribution
exists and is given by $t\mapsto 1_{(t>t_n)}q_{y_n}(t,t-t_n)$, the Radon-Nikodym
derivative of the hazard measure. This proves the result.
\end{namedproof}

\begin{lemma}
\label{lemma:ConditionalJumpDistSMP}
The conditional distribution of $Y_{n+1}$ given $T_0,Y_0,\ldots,Y_n,T_{n+1}$ is 
almost surely given by
\begin{align}
  P(Y_n=j|T_0,Y_0,\ldots,Y_n,T_{n+1})
  &=\frac{q_{Y_nj}(t,t-T_{n+1})}{q_{Y_n}(t,t-T_{n+1})},
\end{align}
understanding that on the almost sure set when the above holds, the denominator is nonzero.
\end{lemma}
\begin{namedproof}{\ref{lemma:ConditionalJumpDistSMP}}
With $\PPP(E)$ denoting the power set of $E$, let $\pi^{(n)}_{z_n,t}:\PPP(E)\to[0,1]$ denote
the conditional distribution of $Y_{n+1}$ given
$(T_0,Y_0,\ldots,T_n,Y_n,T_{n+1})=(z_n,t)$. Let $\Lambda:D(E)\to D(E)$ denote
the total compensator of $Z$, and let the
compensator for the jumps to state $i$ be denoted by $\Lambda^i:D(E)\to D(E)$.
By (4.48) of \cite{MR2189574}, we have
that whenever $q_{y_n}(t,t-t_n)\neq 0$, it holds for all $A\subseteq E$ that
\begin{align}
	\pi^{(n)}_{z_n,t}(A) &= \sum_{i\in A}\frac{\dv \Lambda^i}{\dv \Lambda}(t)
	=\sum_{i\in A} \frac{1_{(y_n\neq i)}q_{y_ni}(t,t-t_n)}{q_{y_n}(t,t-t_n)},
\end{align}
which yields the result.
\end{namedproof}

\begin{lemma}
\label{lemma:ConditionalEventWithTime}
The conditional distribution of $T_{n+1}$ given $Z_t = i$, $U_t = u$ and $N_t = n$
has hazard function $h:\RR_+\to\RR_+$ given by $h(s) = 1_{(s>t)}q_i(t,u+s-t)$.
\end{lemma}
\begin{namedproof}{Lemma \ref{lemma:ConditionalEventWithTime}}
Let $s\ge0$. We find that
\begin{align}
	&P(T_{n+1}\le s|Z_t = i,U_t = u,N_t=n)\notag\\
	&=P(T_{n+1}\le s|Y_n = i,T_n = t-u,T_{n+1}> t)\notag\\
	&=\frac{P(T_{n+1}\le s,T_{n+1}>t|Y_n=i,T_n=t-u)}{P(T_{n+1}> t|Y_n=i,T_n=t-u)}.\label{eq:CondAtTimeHazard1}
\end{align}
Now, by Lemma \ref{lemma:ConditionalEventHazardSMP}, the distribution
of $T_{n+1}$ given $Y_n = i$ and $T_n=t-u$ has hazard function
$v\mapsto 1_{(v>t-u)}q_{i}(v,v-(t-u))$. Therefore, we in particular obtain
\begin{align}
  P(T_{n+1}> t|Y_n=i,T_n=t-u)
  &=\exp\left(\int_{t-u}^t q_i(v,v-(t-u))\dv v\right),\label{eq:CondAtTimeHazard2}
\end{align}
and, for $s>t$,
\begin{align}
  &\frac{\dv}{\dv s}P(T_{n+1}\le s,T_{n+1}> t|Y_n=i,T_n=t-u)\notag\\
  &=q_i(s,s-(t-u))\exp\left(\int_{t-u}^s q_i(v,v-(t-u))\dv v\right).\label{eq:CondAtTimeHazard3}
\end{align}
Combining (\ref{eq:CondAtTimeHazard1}), (\ref{eq:CondAtTimeHazard2}) and (\ref{eq:CondAtTimeHazard3})
yields that for $s>t$, it holds that
\begin{align}
	&P(T_{n+1}\le s|Z_t = i,U_t = u,N_t=n)\notag\\
	&=q_i(s,s-(t-u))\exp\left(\int_{t}^s q_i(v,v-(t-u))\dv v\right),
\end{align}
proving that the conditional distribution of $T_{n+1}$ given $Z_t = i$, $U_t = u$ and $N_t = n$
has hazard $s\mapsto q_i(t,u+s-t)$ for $s>t$. As the support of the distribution
is contained in $[t,\infty)$, this proves the lemma.
\end{namedproof}

\begin{lemma}
\label{lemma:ConditionalNextEventWithTime}
The conditional distribution of $T(t)$ given $Z_t = i$ and $U_t = u$
has hazard function $h:\RR_+\to\RR_+$ given by $h(s) = 1_{(s>t)}q_i(t,u+s-t)$.
\end{lemma}
\begin{namedproof}{Lemma \ref{lemma:ConditionalNextEventWithTime}}
Define $g:\NN_0\to[0,1]$ by $g(n)=P(N_t = n|Z_t = i,U_t = u)$. Applying
Lemma \ref{lemma:ConditionalEventWithTime}, it then holds for $s>t$ that
\begin{align}
	P(T(t)\le s|Z_t = i,U_t = u)
  &=\sum_{n=0}^\infty P(T(t)> s|Z_t = i,U_t = u,N_t = n)g(n)\notag\\
  &=\sum_{n=0}^\infty P(T_{n+1}> s|Z_t = i,U_t = u,N_t = n)g(n)\notag\\
  &=\sum_{n=0}^\infty \exp\left(\int_t^s q_i(v,u+v-t)\dv v\right) g(n)\notag\\
  &=\exp\left(\int_t^s q_i(v,u+v-t)\dv v\right),
\end{align}
yielding the result.
\end{namedproof}

From a heuristic perspective, Lemma \ref{lemma:ConditionalNextEventWithTime} allows
us some insight into the distribution of $Z$. The lemma states that given the values
of $(Z_t,U_t)$, the distribution next jump has a fixed hazard, independently
of the number of jumps made. We are now ready to prove the results stated in
Section \ref{sec:Intensities}.

\begin{namedproof}{Lemma \ref{lemma:SMGTwoJumpsBound}}
For $n\ge0$, define $g(n) = P(N_t = n,Z_t = i, U_t=u)$. We then have
\begin{align}
  &P(N_{t+h}-N_t\ge 2|Z_t = i,U_t=u)\notag\\
  &=\sum_{n=0}^\infty P(N_{t+h}-N_t\ge 2|Z_t = i,U_t=u,N_t=n)g(n)\notag\\
  &=\sum_{n=0}^\infty P(t<T_{n+1}<T_{n+2}\le t+h|Z_t = i,U_t=u,N_t=n)g(n).\label{eq:SMGTwoJumpsBoundMain}
\end{align}
Now, by Lemma \ref{lemma:ConditionalEventWithTime}, the conditional distribution
of $T_{n+1}$ given $Z_t=i$, $U_t=u$ and $N_t = n$ has hazard function
$q_i(s,u+s-t)$ for $s\ge t$. Applying Lemma \ref{lemma:ConditionalEventHazardSMP},
we then obtain that the conditional distribution of $Z_{T_{n+1}}$ and $T_{n+1}$
given $Z_t = i$, $U_t = u$ and $N_t = n$ has density
\begin{align}
  f(k,s) &= q_{ik}(s,u+s-t)\exp\left(-\int_t^s q_i(v,u+v-t)\dv v\right).
\end{align}
for $s\ge t$ and zero otherwise. Next, noting that we always have $U_{T_{n+1}}=0$,
we may use the strong Markov property of $(Z,U)$ at $T_{n+1}$ and obtain
\begin{align}
  &P(t<T_{n+1}<T_{n+2}\le t+h|Z_t = i,U_t=u,N_t=n)\notag\\
  &=\sum_{k\in E}\int_t^{t+h} P(s<T_{n+2} \le t+h|Z_{T_{n+1}}=k,T_{n+1}=s)f(k,s)\dv s\label{eq:TwoJumpsCondComp}.
\end{align}
Here, the distribution of $T_{n+2}$ given $Z_{T_{n+1}}=k$ and $T_{n+1}=s$ has hazard $q_k(v,v-s)$
for $v\ge s$ and zero otherwise, which yields the conditional density
\begin{align}
f(v|k,s) &= q_k(v,v-s)\exp\left(-\int_s^v q_k(r,r-s)\dv r\right).
\end{align}
for $v\ge s$ and zero otherwise. We thus obtain
\begin{align}
  P(s< T_{n+2}\le t+h|Z_{T_{n+1}}=k,T_{n+1}=s)
  &=\int_s^{t+h} f(v|k,s)\dv v\notag\\
  &\le \int_s^{t+h} \|Q(v,v-s)\|_\infty \dv v.
\end{align}
Inserting this in (\ref{eq:TwoJumpsCondComp}), we obtain the bound
\begin{align}
  &P(t<T_{n+1}<T_{n+2}\le t+h|Z_t = i,U_t=u,N_t=n)\notag\\
  &\le \sum_{k\in E}\int_t^{t+h}\int_s^{t+h} \|Q(v,v-s)\|_\infty \|Q(s,u+s-t)\|_\infty \dv v \dv s.
\end{align}
As this bound is independent of $n$, we may use it in (\ref{eq:SMGTwoJumpsBoundMain})
and obtain the required result.
\end{namedproof}

\begin{namedproof}{Lemma \ref{lemma:SMGNoTwoJumps}}
By Lemma \ref{lemma:SMGTwoJumpsBound}, we have
\begin{align}
  &\frac{1}{h}P(N_{t+h}-N_t\ge 2|Z_t = i,U_t=u)\notag\\
  &\le h \sup_{t\le s,v\le t+h} \|Q(s,u+s-t)\|_\infty\|Q(v,v-s)\|_\infty.
\end{align}
Now, by the bound (\ref{eq:FiniteMatrixIntensityBound}), we find that $(s,v)\mapsto \|Q(s,v)\|_\infty$
is bounded in both a neighborhood of $(t,u)$ and $(t,0)$. Therefore, the above tends
to zero as $h$ tends to zero.
\end{namedproof}

\begin{namedproof}{Theorem \ref{theorem:IntensityIsTransProbLimitBounded}}
We will show the result with $\eps=1$. To obtain (\ref{eq:IntensityTransProbLimitBound})
for this case, we first fix $t,u\ge0$ and $i\in E$ and seek to bound
the sum of $\frac{1}{h}(p_{ij}(t,t+h,u)-I_{ij})$ over $j\neq i$ for $h\le 1$. Fix $i\neq j$, we then have
\begin{align}
  p_{ij}(t,t+h,u)
  &\le P(N_{t+h} - N_t = 1,Z_{t+h}=j| Z_t = i,U_t = u)\notag\\
   &+ P(N_{t+h} - N_t \ge 2,Z_{t+h}=j| Z_t = i,U_t = u).\label{eq:NoTwoJumpsHelper}
\end{align}
Next, note that by Lemma \ref{lemma:ConditionalJumpDistSMP}
and Lemma \ref{lemma:ConditionalNextEventWithTime}, we have
\begin{align}
  &P(N_{t+h} - N_t = 1,Z_{t+h}=j| Z_t = i,U_t = u)\notag\\
  &\le \int_t^{t+h} q_{ij}(s,u+s-t) \exp\left(-\int_t^s q_i(v,u+v-t)\dv v\right)\dv s\notag\\
  &\le \int_t^{t+h} q_{ij}(s,u+s-t) \dv s,\label{eq:IntensityTransLimitBound1}
\end{align}
since $N_{t+h}-N_t = 1$ and $Z_{t+h}=j$ is equivalent to having the next jump strictly after
$t$ occurring in the interval $(t,t+h]$ with destination state $j$, and the jump
after that occurring strictly after $t+h$. Applying this bound in
(\ref{eq:NoTwoJumpsHelper}), the monotone convergence theorem yields
\begin{align}
\sum_{j\neq i}p_{ij}(t,t+h,u)
&\le \int_t^{t+h} q_i(s,u+s-t) \dv s \notag\\
&+ P(N_{t+h} - N_t \ge 2| Z_t = i,U_t = u).\label{eq:NoTwoJumpsHelper2}
\end{align}
Here, by Lemma \ref{lemma:SMGTwoJumpsBound}, we have
\begin{align}
  &P(N_{t+h}-N_t\ge 2|Z_t = i,U_t=u)\notag\\
  &\le \int_t^{t+h} \int_s^{t+h} \|Q(s,u+s-t)\|_\infty\|Q(v,v-s)\|_\infty \dv v\dv s.\label{eq:IntensityTransLimitBound2}
\end{align}
Combining (\ref{eq:IntensityTransLimitBound1}) and (\ref{eq:IntensityTransLimitBound2})
with (\ref{eq:NoTwoJumpsHelper}), we obtain for $h\le 1$ that
\begin{align}
  &\sum_{j\neq i}\frac{1}{h}p_{ij}(t,t+h,u)\notag\\
  &\le \frac{1}{h}\int_t^{t+h} \|Q(s,u+s-t)\|_\infty\left(1+\int_s^{t+h} \|Q(v,v-s)\|_\infty \dv v\right)\dv s\notag\\
  &\le \sup_{t\le s\le t+1} \|Q(s,u+s-t)\|_\infty\left(1+\sup_{s\le v\le t+1}\|Q(v,v-s)\|_\infty \right)\dv s \label{eq:IntensityTransLimitBound3}.
\end{align}
As we also have
\begin{align}
  \left|\frac{1}{h}p_{ii}(t,t+h,u)-1\right| &= \sum_{j\neq i} \frac{1}{h}p_{ij}(t,t+h,u),
\end{align}
we finally obtain
\begin{align}
  &\sup_{i\in E}\sum_{j\in E}\left|\frac{1}{h}(p_{ij}(t,t+h,u)-I_{ij})\right|\notag\\
  &\le 2\sup_{t\le s\le t+1} \|Q(s,u+s-t)\|_\infty\left(1+\sup_{s\le v\le t+1}\|Q(v,v-s)\|_\infty \right)\dv s\label{eq:IntensityTransLimitBound4}.
\end{align}
Letting $C(t,u)$ be the right-hand side of (\ref{eq:IntensityTransLimitBound3}), we obtain
(\ref{eq:IntensityTransProbLimitBound}) for the case $i\neq j$. By
the right-continuity of the intensities, this definition of the mapping $C$ is measurable,
and by (\ref{eq:FiniteMatrixIntensityBound}), it is bounded on compacts.
\end{namedproof}

\begin{namedproof}{Theorem \ref{theorem:IntensityIsTransProbLimit}}
Consider first the case $i\neq j$. Here, we need to show that
\begin{align}
\label{eq:IntensityTransProbLimitOffDiag}
  q_{ij}(t,u) &= \lim_{h\to0}\frac{1}{h}p_{ij}(t,t+h,u),
\end{align}
To this end, we first make the decomposition
\begin{align}
  p_{ij}(t,t+h,u)
  &=P(N_{t+h} - N_t = 1, Z_{t+h} = j | Z_t = i,U_t = u)\notag\\
  &+P(N_{t+h} - N_t \ge 2, Z_{t+h} = j | Z_t = i,U_t = u).\label{eq:IntensityIsTransLimitHelper}
\end{align}
Here, Lemma \ref{lemma:SMGNoTwoJumps} shows that the second term tends
to zero, so it suffices to show that the first term tends to $q_{ij}(t,u)$. To
do so, we note that having $N_{t+h}-N_t=1$ and $Z_{t+h}=j$ is equivalent to
having $T(t)\le t+h$, $T(T(t))>t+h$ and $Z_{T(t)}=j$. Defining
\begin{align}
	f(t,u,s) &= q_{ij}(s,u+s-t) \exp\left(-\int_t^s q_i(v,u+v-t)\dv v\right),\\
	a(t,r) &= \exp\left(-\int_t^r q_j(v,v)\dv v\right),
\end{align}
for $s\ge t$ and $r\ge t$, we find by Lemma \ref{lemma:ConditionalNextEventWithTime} and
Lemma \ref{lemma:ConditionalJumpDistSMP} that $s\mapsto f(t,u,s)$ is the
conditional density given $Z_t = i$ and $U_t = u$ of making the next jump at time $t$
to state $j$, and $r\mapsto a(t,r)$ is the survival function for the next
jump given $Z_t = j$ and $U_t = 0$. Therefore, we obtain
\begin{align}
  &=P(N_{t+h} - N_t = 1, Z_{t+h} = j | Z_t = i,U_t = u)\notag\\
  &=\int_t^{t+h} f(t,u,s)a(s,t+h) \dv s.
\end{align}
As the intensities are assumed to be right-continuous, this yields
\begin{align}
  \lim_{h\to0}\frac{1}{h} P(N_{t+h} - N_t = 1, Z_{t+h} = j | Z_t = i,U_t = u)
  &= q_{ij}(t,u),
\end{align}
which, combined with (\ref{eq:IntensityIsTransLimitHelper}), yields
(\ref{eq:IntensityTransProbLimitOffDiag}). It remains to prove for any $i\in E$ that
\begin{align}
\label{eq:IntensityTransProbLimitDiag}
  q_{ii}(t,u) &= \lim_{h\to0}\frac{1}{h}(p_{ii}(t,t+h,u)-1).
\end{align}
In order to obtain this, we simply apply the dominated convergence theorem
with the bound $C(t,u)$ from Theorem \ref{theorem:IntensityIsTransProbLimitBounded}
to obtain
\begin{align}
	\lim_{h\to0}\frac{1}{h}(p_{ii}(t,t+h,u)-1)
	&=-\lim_{h\to0}\sum_{j\neq i} \frac{1}{h} p_{ij}(t,t+h,u)\notag\\
	&=-\sum_{j\neq i} \lim_{h\to0} \frac{1}{h} p_{ij}(t,t+h,u)\notag\\
	&=-\sum_{j\neq i} q_{ij}(t,u) = q_{ii}(t,u),
\end{align}
as required.
\end{namedproof}

\subsection{Proofs for Section \ref{sec:Forward}}

We begin by proving a few auxiliary lemmas, including a version of the
Chapman-Kolmogorov equations for semi-Markov processes.

\begin{lemma}
\label{lemma:DurationPartiallyNonExpansive}
It holds for all $0\le s\le t$ that $U_t - U_s\le t-s$.
\end{lemma}
\begin{namedproof}{Lemma \ref{lemma:DurationPartiallyNonExpansive}}
Fix $0\le s\le t$. If there are no event times in $(s,t]$, then Lemma \ref{lemma:SMGSimplePathProperties}
shows that $U_t-U_s = t-s$. Considering the case where there is an event time in $(s,t]$, let
$T$ be the largest such event time, $T=\sup\{T_n| s<T_n\le t\}$. We then obtain
\begin{align}
	U_t - U_s &= t-T - U_s \le t-T\le t-s,
\end{align}
as required, since $U$ is nonnegative.
\end{namedproof}
		
\begin{lemma}
\label{lemma:pijSMGMeasureSupport}
Let $0\le s\le t$, and let $u\ge0$. For $i\neq j$, it holds that the measure $p_{ij}(s,t,u,\cdot)$ is concentrated
on the set $[0,t-s]$. For $i\in E$, it holds that the measure $p_{ii}(s,t,u,\cdot)$ is concentrated
on $[0,t-s]\cup \{u,u+t-s\}$.
\end{lemma}

\begin{namedproof}{Lemma \ref{lemma:pijSMGMeasureSupport}}
Consider first the case where $i\neq j$. For $s=t$, it is immediate that $p_{ij}(s,t,u,\cdot)$
is the zero measure, so the result is immediate in this case. Consider instead
the case where $s<t$. Here, we have
\begin{align}
p_{ij}(s,t,u,A) &= P(Z_t = j,U_t\in A|Z_s = i,U_s = u)\notag\\
&= P(Z_s = i, Z_t = j,U_t\in A|Z_s = i,U_s = u).\label{eq:pijSMGUpperBoundHelper}
\end{align}
Now, when $Z_s = i$ and $Z_t = j$, with $s<t$, it must hold that there is some event
time for $Z$ in the interval $(s,t]$. Therefore, we may apply
Lemma \ref{lemma:DurationPartiallyNonExpansive} and obtain
\begin{align}
	(Z_s = i,Z_t = j)
	&\subseteq \cup_{n=1}^\infty (s<T_n\le t,Z_s = i,Z_t = j)\notag\\
	&\subseteq \cup_{n=1}^\infty (s<T_n\le t,U_{T_n} = 0, Z_s = i,Z_t = j)\notag\\
	&\subseteq \cup_{n=1}^\infty (s<T_n\le t,U_t \le t-T_n,U_{T_n} = 0, Z_s = i,Z_t = j)\notag\\
	&\subseteq (U_t \le t-s).\label{eq:UpperBoundHelper2}
\end{align}
Using this in (\ref{eq:pijSMGUpperBoundHelper}), we obtain
\begin{align}
	p_{ij}(s,t,u,A) &= P(Z_t = j,U_t\in A|Z_s = i,U_s = u)\notag\\
	&= P(Z_s = i, Z_t = j,U_t\in A\cap[0,t-s]|Z_s = i,U_s = u),
\end{align}
as required. Next, consider a single $i\in E$. In the case where $s=t$, it is immediate
that $p_{ii}(s,t,u,\cdot)$ is concentrated on $\{u\}$, and so the result follows. Consider
instead $s<t$. In this case, it either holds that there are zero event times in $(s,t]$,
or there is an event time. In the latter case, we can use a calculation as in
(\ref{eq:UpperBoundHelper2}) to obtain that the duration at time $t$ must be in
$[0,u+t-s]$. In the former case, we have
\begin{align}
	&\cap_{n=1}^\infty (T_n\notin(s,t], Z_s = i,U_s = u, Z_t = j)\notag\\
  &=\cap_{n=1}^\infty (T_n\notin(s,t], Z_s = i,U_s = u, Z_t = j,U_t = u+t-s).
\end{align}
From this, the result follows.
\end{namedproof}

The analogue of the Chapman-Kolmogorov equations for semi-Markov processes is encapsulated
in the following lemma. Note that the integration domain in (\ref{eq:SMGChapmanKolmo}) is
assumed to be $[0,u+t-s]$ and not $[0,u+t-s)$. This is important as $p_{ij}(s,r,u,\cdot)$
in the case $i=j$ can have a point mass in $u+t-s$, corresponding to the case where no
jumps are made in the interval $(s,r]$.

\begin{lemma}
\label{lemma:SMGChapmanKolmo}
Fix $0\le s\le r\le t$. It holds that
\begin{align}
\label{eq:SMGChapmanKolmo}
p_{ij}(s,t,u,A)&=\sum_{k\in E}\int_0^{u+r-s} p_{kj}(r,t,v,A) p_{ik}(s,r,u,\dv v).
\end{align}
In particular, we have for $d\ge0$ that
\begin{align}
\label{eq:SMGChapmanKolmoCDF}
p_{ij}(s,t,u,d)&=\sum_{k\in E}\int_0^{u+r-s} p_{kj}(r,t,v,d) p_{ik}(s,r,u,\dv v).
\end{align}
\end{lemma}

\begin{namedproof}{Lemma \ref{lemma:SMGChapmanKolmo}}
By the Chapman-Kolmogorov equations for the inhomogeneous Markov chain $(Z,U)$, we have
\begin{align}
p_{ij}(s,t,u,A)
&=P_{s,t}(i,u,\{j\}\times A)
=\int_{E\times\RR_+} P_{r,t}(k,v,\{j\}\times A)\dv P_{s,r}(i,u,\dv k,\dv v)\notag\\
&=\sum_{k\in E}\int_0^\infty P_{r,t}(k,v,\{j\}\times A)\dv P_{s,r}(i,u,k,\dv v)\notag\\
&=\sum_{k\in E}\int_0^\infty p_{kj}(r,t,v,A) p_{ik}(s,r,u,\dv v).
\end{align}
As Lemma \ref{lemma:pijSMGMeasureSupport} shows that the support of $p_{ik}(s,r,u,\cdot)$
is included in $[0,u+r-s]$, we obtain the result. The identity (\ref{eq:SMGChapmanKolmoCDF})
follows from (\ref{eq:SMGChapmanKolmo}) by definition.
\end{namedproof}

\begin{namedproof}{Lemma \ref{lemma:SMGNoQuickCycles}}
We first note that for $Z$ to remain in state $i$ at time $t+h$ conditionally on being
in state $i$ at time $t$ can only happen by changing state an even number of times. In particular,
this yields
\begin{align}
  p_{ii}(t,t+h,u,(u+h)-)
  &=P(Z_{t+h} = i,U_{t+h} < u+h|Z_t = i,U_t = u) \\
  &\le P(N_{t+h}-N_t = 0,U_{t+h} < u+h|Z_t = i,U_t = u)\notag\\
  &+P(N_{t+h}-N_t \ge 2|Z_t = i,U_t = u).\notag
\end{align}
Here, we have
\begin{align}
  &P(N_{t+h}-N_t = 0,U_{t+h}<u+h|Z_t = i,U_t = u)\notag\\
  &=P(N_{t+h}-N_t = 0,U_t=u, U_{t+h} < u+h|Z_t = i,U_t = u)=0,
\end{align}
since, when making no jumps in the time interval $(t,t+h]$, $U_t = u$
implies that $U_{t+h}=u+h$. As a consequence, we obtain
\begin{align}
  \frac{1}{h} p_{ii}(t,t+h,u,(u+h)-)
  &\le \frac{1}{h} P(N_{t+h}-N_t \ge 2|Z_t = i,U_t = u).
\end{align}
The limit statement (\ref{eq:NoQuickCyclesLimit}) then follows from Lemma \ref{lemma:SMGNoTwoJumps},
and the existence of the bound (\ref{eq:NoQuickCyclesBound}) follows as in
the proof of Theorem \ref{theorem:IntensityIsTransProbLimitBounded}.
\end{namedproof}

In order to prove Lemma \ref{lemma:pIsDurationDiff}, we first define
\begin{align}
	r_{ij}(s,u,t)&=q_{ij}(t,u+t-s)\exp\left(-\int_s^t q_i(v,u+v-s)\dv v\right),\label{eq:rij}\\
	R_i(s,u,t)&=\exp\left(-\int_s^t q_i(v,u+v-s)\dv v\right).\label{eq:Ri}
\end{align}
Note that $(j,t)\mapsto r_{ij}(s,u,t)$ is then the conditional density of the
next jump and its destination state given $Z_s=i$ and $U_s = u$, and $t\mapsto R_i(s,u,t)$
is the survival function for the next jump given $Z_s = i$ and $U_s = u$.

\begin{namedproof}{Lemma \ref{lemma:pIsDurationDiff}}
First consider the case where $i\neq j$. Note that as
$p_{ij}(s,t,u,\cdot)$ is concentrated on $[0,t-s]$ by Lemma \ref{lemma:pijSMGMeasureSupport}, we have that whenever
$d-h\ge t-s$, it holds that
\begin{align}
  p_{ij}(s,t,u,d) - p_{ij}(s,t,u,d-h)
  &=p(s,t,u,(d-h,d])=0.
\end{align}
In particular, if $d>t-s$, it is immediate that the derivative from the left exists
and is zero. Next, consider the case where $d\le t-s$. As $s\le t-d$ in this case,
we can apply the Chapman-Kolmogorov equations to obtain
\begin{align}
  &\lim_{h\to0}\frac{1}{h}p_{ij}(s,t,u,d) - p_{ij}(s,t,u,d-h)\notag\\
  &=\lim_{h\to0}\frac{1}{h}p(s,t,u,(d-h,d])\notag\\
  &=\sum_{k\in E}\lim_{h\to0}\frac{1}{h} \int_0^{u+t-d-s}p_{kj}(t-d,t,v,(d-h,d]) p_{ik}(s,t-d,u,\dv v),\label{eq:DurationCKHelp}
\end{align}
insofar as the limits exist. We wish to argue that for each $k$ in the sum above, the limit exists.
To this end, we first consider the case $k=j$. Note that if $U_{t-d} = v$ for $v>0$ and there are no jumps on $(t-d,t]$, then
$U_t >d$. Therefore, we obtain
\begin{align}
  &p_{jj}(t-d,t,v,(d-h,d])\notag\\
  &=P(Z_t = j,d-h<U_t\le d| Z_{t-d} = j,U_{t-d} = v)\notag\\
  &= P(N_t-N_{t-d} \ge 2, Z_t = j,d-h<U_t\le d| Z_{t-d} = j,U_{t-d} = v).
\end{align}
Here, Lemma \ref{lemma:SMGNoTwoJumps} yields
\begin{align}
  &\limsup_{h\to0}\frac{1}{h}P(N_t-N_{t-d} \ge 2, Z_t = j,d-h<U_t\le d| Z_{t-d} = j,U_{t-d} = v)\notag\\
  &\le \limsup_{h\to0}\frac{1}{h}P(N_t-N_{t-d} \ge 2| Z_{t-d} = j,U_{t-d} = v)=0.
\end{align}
Therefore, applying Lemma \ref{lemma:SMGTwoJumpsBound} and the dominated convergence
theorem, we obtain
\begin{align}
  \lim_{h\to0}\frac{1}{h} \int_0^{u+t-d-s}p_{jj}(t-d,t,v,(d-h,d]) p_{ij}(s,t-d,u,\dv v)
  =0.\label{eq:DurationCKDiagHelp}
\end{align}
Next, consider the case $k\neq j$. In this case, Lemma \ref{lemma:pijSMGMeasureSupport} yields that $p_{kj}(t-d,t,v,\cdot)$ is concentrated
on $[0,d]$. Also note that when $U_t>d-h$, no jumps are made in the time interval $(t-d+h,t]$.
Therefore, any jumps made by the process in the interval $(t-d,t]$ must in this case be made in the
interval $(t-d,t-d+h]$, yielding
\begin{align}
  &p_{kj}(t-d,t,v,(d-h,d])\notag\\
  &=P(Z_t = j, d-h<U_t|Z_{t-d} = k, U_{t-d} = v)\notag\\
  &=P(N_{t-d+h}-N_{t-d} =1,Z_t = j, d-h<U_t|Z_{t-d} = k, U_{t-d} = v)\notag\\
  &+P(N_{t-d+h}-N_{t-d} \ge 2,Z_t = j, d-h<U_t|Z_{t-d} = k, U_{t-d} = v)\notag.
\end{align}
As in the previous case, Lemma \ref{lemma:SMGNoTwoJumps} yields
\begin{align}
  &\limsup_{h\to0}\frac{1}{h}P(N_{t-d+h}-N_{t-d} \ge 2,Z_t = j, d-h<U_t|Z_{t-d} = k, U_{t-d} = v)\notag\\
  &\le \limsup_{h\to0}\frac{1}{h}P(N_{t-d+h}-N_{t-d} \ge 2|Z_{t-d} = k, U_{t-d} = v)=0,
\end{align}
with the convergence, as before, being bounded above by a constant according to 
Lemma \ref{lemma:SMGTwoJumpsBound}. Next, note that on $N_{t-d} = n$,
having $N_{t-d+h}-N_{t-d} =1$ and $d-h<U_t$ is equal to having
$t-d< T_{n+1}\le t-d+h$ and $T_{n+2}>t$.  Let $g(n) = P(N_{t-d} = n|Z_{t-d}=k,U_{t-d}=v)$
and $A = (t-d,t-d+h]$, we therefore obtain
\begin{align}
  &P(N_{t-d+h}-N_{t-d} =1,Z_t = j, d-h<U_t|Z_{t-d} = k, U_{t-d} = v)\notag\\
  &=\sum_{n=0}^\infty P(T_{n+1}\in A,T_{n+2}>t,Z_t = j|Z_{t-d} = k, U_{t-d} = v,N_{t-d} = n)g(n).\label{eq:durationDerivHelper}
\end{align}
Here, by Lemma \ref{lemma:ConditionalEventWithTime} and Lemma \ref{lemma:ConditionalJumpDistSMP},
the conditional distribution of $(T_{n+1},Z_{T_{n+1}})$ given $Z_{t-d} = k$
$U_{t-d} = v$ and $N_{t-d} = n$ has density $(l,s)\mapsto r_{kl}(t-d,v,s)$ for $l\neq k$
and $s\ge t-d$, and by Lemma \ref{lemma:ConditionalEventHazardSMP}, the survival function
of the conditional distribution of $T_{n+2}$ given $T_{n+1}=s$ and $Z_{T_{n+1}}=l$ is $w\mapsto R_l(s,0,w)$.
Here, we use the notation outlined in (\ref{eq:rij}) and (\ref{eq:Ri}). Thus, we obtain
\begin{align}
  &P(T_{n+1}\in A,T_{n+2}>t,Z_t = j|Z_{t-d} = k, U_{t-d} = v,N_{t-d} = n)\notag\\
  &=\int_{t-d}^{t-d+h} r_{kj}(t-d,v,s)R_j(s,0,t)\dv s.
\end{align}
As this is independent of $n$, insertion in (\ref{eq:durationDerivHelper}) yields
\begin{align}
  &P(N_{t-d+h}-N_{t-d} =1,Z_t = j, d-h<U_t|Z_{t-d} = k, U_{t-d} = v)\notag\\
  &=\int_{t-d}^{t-d+h} R_j(s,0,t)r_{kj}(t-d,v,s)\dv s.\label{eq:durDiffMainTerm}
\end{align}
By right-continuity of the intensities, this yields
\begin{align}
  &\lim_{h\to0}\frac{1}{h}P(N_{t-d+h}-N_{t-d} =1,Z_t = j, d-h<U_t|Z_{t-d} = k, U_{t-d} = v)\notag\\
  &=r_{kj}(t-d,v,t-d)R_j(t-d,0,t)\notag\\
  &=q_{kj}(t-d,v)\exp\left(-\int_{t-d}^t q_j(r,r-(t-d))\dv r\right).\label{eq:DurationCKOffDiagHelp}
\end{align}
For brevity, we let $a_{kj}(t,d,u,v)$ denote the right-hand side of
(\ref{eq:DurationCKOffDiagHelp}). Now using (\ref{eq:DurationCKDiagHelp}) and (\ref{eq:DurationCKOffDiagHelp})
in (\ref{eq:DurationCKHelp}), we obtain
\begin{align}
  &\lim_{h\to0}\frac{1}{h}p_{ij}(s,t,u,d) - p_{ij}(s,t,u,d-h)\notag\\
  &=\sum_{k\neq j}\lim_{h\to0}\frac{1}{h} \int_0^{u+t-d-s}p_{kj}(t-d,t,v,(d-h,d]) p_{ik}(s,t-d,u,\dv v)\notag\\
  &=\sum_{k\neq j} \int_0^{u+t-d-s} a_{kj}(t,d,u,v) p_{ik}(s,t-d,u,\dv v),\label{eq:DurDiffFinal}
\end{align}
yielding (\ref{eq:TransProbDurationDiffNonDiag}). Note that in order to move the limit
under the integral in (\ref{eq:DurDiffFinal}), we applied the dominated convergence
theorem, making use of that the integral in (\ref{eq:durDiffMainTerm}) is bounded
from above on $[t-d,t-d+\eps]$ for some $\eps>0$. This shows that $p_{ij}(s,t,u,d)$ is differentiable
from the left in $d$ for $d>0$ when $i\neq j$. The diagonal case follows immediately from this
by writing $p_{ii}(s,t,u,d)$ as one minus the sum of $p_{ij}(s,t,u,d)$ for $j\neq i$.
\end{namedproof}

\begin{namedproof}{Theorem \ref{theorem:SMPForward}}
Let $s,u\ge0$ and let $t\ge s$ and $d>0$. Applying the Chapman-Kolmogorov equations of Lemma \ref{lemma:SMGChapmanKolmo}, we have
\begin{align}
&\frac{\partial p_{ij}}{\partial t}(s,t,u,d)\notag\\
&=\lim_{h\to0}\frac{1}{h}(p_{ij}(s,t+h,u,d)-p_{ij}(s,t,u,d))\notag\\
&=\lim_{h\to0}\frac{1}{h}\sum_{k\in E}\int_0^{u+t-s} p_{kj}(t,t+h,v,d) p_{ik}(s,t,u,\dv v)-\frac{1}{h}p_{ij}(s,t,u,d),
\end{align}
insofar as the limits exist, our objective is to argue that this is the case.
Now, for $k\neq j$, Lemma \ref{lemma:pijSMGMeasureSupport} shows that
the measure $p_{kj}(t,t+h,v,\cdot)$ is concentrated on $[0,h]$. Therefore,
whenever $d\ge h$, we have $p_{kj}(t,t+h,v,d) = p_{kj}(t,t+h,v)$. Applying
Lemma \ref{theorem:IntensityIsTransProbLimit}, Lemma \ref{theorem:IntensityIsTransProbLimitBounded} and
the dominated convergence theorem, we then obtain
\begin{align}
&\lim_{h\to0} \int_0^{u+t-s} \frac{1}{h} p_{kj}(t,t+h,v,d) p_{ik}(s,t,u,\dv v)\notag\\
&=\int_0^{u+t-s} \lim_{h\to0} \frac{1}{h} p_{kj}(t,t+h,v) p_{ik}(s,t,u,\dv v)\notag\\
&=\int_0^{u+t-s} q_{kj}(t,v) p_{ik}(s,t,u,\dv v).\label{eq:SMGForwardOffDiagonal}
\end{align}
Now consider the case $k=j$. We wish to evaluate the limit
\begin{align}
\label{eq:SMGForwardDiagonalLimit}
	\lim_{h\to0} \int_0^{u+t-s} \frac{1}{h} p_{jj}(t,t+h,v,d) p_{ij}(s,t,u,\dv v)-\frac{1}{h}p_{ij}(s,t,u,d).
\end{align}
We note that for $v\ge d$, Lemma \ref{lemma:SMGNoQuickCycles} yields
\begin{align}
  \limsup_{h\to0} \frac{1}{h} p_{jj}(t,t+h,v,d)
  &\le \limsup_{h\to0} \frac{1}{h} p_{jj}(t,t+h,v,(v+h)-) = 0,
\end{align}
so that $p_{jj}(t,t+h,v,d)/h$ tends to zero as $h$ tends to zero from above. The bound in Lemma
\ref{lemma:SMGNoQuickCycles} and the dominated convergence theorem then allows us to conclude
that in the case $d< u+t-s$, we have
\begin{align}
	&\lim_{h\to0} \int_d^{u+t-s} \frac{1}{h} p_{jj}(t,t+h,v,d) p_{ij}(s,t,u,\dv v)=0.
\end{align}
Therefore, we obtain that if $d< u+t-s$, then the limit (\ref{eq:SMGForwardDiagonalLimit}) is equal
to the limit
\begin{align}
\label{eq:SMGForwardDiagonalLimit2}
	\lim_{h\to0} \int_0^d \frac{1}{h} p_{jj}(t,t+h,v,d) p_{ij}(s,t,u,\dv v)-\frac{1}{h}p_{ij}(s,t,u,d),
\end{align}
provided that this exists. In the case $d\ge u+t-s$, we we note that according to Lemma \ref{lemma:pijSMGMeasureSupport},
$p_{ij}(s,t,u,\cdot)$ is concentrated on $[0,u+t-s]$, so the limits
(\ref{eq:SMGForwardDiagonalLimit}) and (\ref{eq:SMGForwardDiagonalLimit2}) are also
equal in this case, whenever they exist.

We proceed with evaluating (\ref{eq:SMGForwardDiagonalLimit2}). Note that
$p_{jj}(t,t+h,v,\cdot)$ is concentrated on $[0,v+h]$ by Lemma \ref{lemma:pijSMGMeasureSupport}.
Therefore, $p_{jj}(t,t+h,v,d)=p_{jj}(t,t+h,v)$ for $0\le v \le d-h$, for $d-h< v\le d$, we have
\begin{align}
	p_{jj}(t,t+h,v,d)
    &=p_{jj}(t,t+h,v,v+h)-p_{jj}(t,t+h,v,(d,v+h])\notag\\
    &=p_{jj}(t,t+h,v)-p_{jj}(t,t+h,v,(d,v+h]).
\end{align}
As a consequence, we find that
\begin{align}
	&\int_0^d \frac{1}{h} p_{jj}(t,t+h,v,d) p_{ij}(s,t,u,\dv v)-\frac{1}{h}p_{ij}(s,t,u,d)\notag\\
  &=\int_0^d \frac{1}{h} p_{jj}(t,t+h,v) p_{ij}(s,t,u,\dv v)-\frac{1}{h}p_{ij}(s,t,u,d)\notag\\
  &-\frac{1}{h} \int_{d-h}^d p_{jj}(t,t+h,v,(d,v+h]) p_{ij}(s,t,u,\dv v),\label{eq:SMGForwardDiagonalDecomp}
\end{align}
where the final integral is over $(d-h,d]$. Taking the limit of the first two terms in (\ref{eq:SMGForwardDiagonalDecomp}), we may
apply the dominated convergence theorem as in the off-diagonal case and obtain
\begin{align}
&\lim_{h\to0} \int_0^d \frac{1}{h} p_{jj}(t,t+h,v) p_{ij}(s,t,u,\dv v)-\frac{1}{h}p_{ij}(s,t,u,d)\notag\\
&=\lim_{h\to0} \int_0^d \frac{1}{h} (p_{jj}(t,t+h,v)-1) p_{ij}(s,t,u,\dv v)\notag\\
&= \int_0^d q_{jj}(t,v) p_{ij}(s,t,u,\dv v).\label{eq:SMGForwardSimpleDiagonalPart}
\end{align}
As regards the term being subtracted in (\ref{eq:SMGForwardDiagonalDecomp}), we can write
\begin{align}
	&\frac{1}{h} \int_{d-h}^d p_{jj}(t,t+h,v,(d,v+h]) p_{ij}(s,t,u,\dv v)\notag\\
  &=\frac{1}{h}(p_{ij}(s,t,u,d) - p_{ij}(s,t,u,d-h))\notag\\
  &+\frac{1}{h}\int_{d-h}^d (p_{jj}(t,t+h,v,(d,v+h])-1) p_{ij}(s,t,u,\dv v).\label{eq:SMGForwardFurtherDecomp}
\end{align}
For the first term in (\ref{eq:SMGForwardFurtherDecomp}), we may
use Lemma \ref{lemma:pIsDurationDiff} to obtain that the limit as $h$ tends to zero exists
and is given by
\begin{align}
\label{eq:SMGForwardPartialTerm}
  \lim_{h\to}\frac{1}{h}(p_{ij}(s,t,u,d) - p_{ij}(s,t,u,d-h))
  &=\frac{\partial p_{ij}}{\partial d}(s,t,u,d),
\end{align}
where the right-hand side denotes the left derivative. As for the second term in (\ref{eq:SMGForwardFurtherDecomp}), we have
\begin{align}
  &\left|\frac{1}{h} \int_{d-h}^d (p_{jj}(t,t+h,v,(d,v+h])-1) p_{ij}(s,t,u,\dv v)\right|\notag\\
  &=\left|\int_{d-h}^d \frac{1}{h} (p_{jj}(t,t+h,v)-1)+\frac{1}{h}p_{jj}(t,t+h,v,d) p_{ij}(s,t,u,\dv v)\right|\notag\\
  &\le \left|\int_{d-h}^d \frac{1}{h} (p_{jj}(t,t+h,v)-1)p_{ij}(s,t,u,\dv v)\right|\notag\\
  &+\left|\int_{d-h}^d \frac{1}{h} p_{jj}(t,t+h,v,d) p_{ij}(s,t,u,\dv v)\right|.
\end{align}
Here, the integrand in the first term of the final estimate converges to $q_{jj}(t,v)$,
and by Theorem \ref{theorem:IntensityIsTransProbLimitBounded}, we also have for
$d-h<v\le d$ that
\begin{align}
\frac{1}{h} (p_{jj}(t,t+h,v)-1)
  &\le \sup_{d-h<v\le d}C(t,u),
\end{align}
the latter being finite. Therefore, the dominated convergence theorem yields that the first
term converges to zero. As regards the second integrand, we note that $d-h< v\le d$ implies $d<v+h$, so
Lemma \ref{lemma:SMGNoQuickCycles} yields that the integrand in this case tends to zero,
dominated by a bound measurable in $t$ and $u$ and bounded on compacts. Therefore,
the dominated convergence theorem also allows us to conclude that the second integral
also tends to zero. All in all, this yields
\begin{align}
\label{eq:SMGForwardFinalRemainder}
 \lim_{h\to0}\left|\int_{d-h}^d \frac{1}{h} (p_{jj}(t,t+h,v,(d,v+h])-1) p_{ij}(s,t,u,\dv v)\right|=0.
\end{align}
Combining (\ref{eq:SMGForwardFinalRemainder}) with our previous conclusions
and simplifications, we finally conclude
\begin{align}
&\lim_{h\to0} \int_0^{u+t-s} \frac{1}{h} p_{jj}(t,t+h,v,d) p_{ij}(s,t,u,\dv v)-\frac{1}{h}p_{ij}(s,t,u,d)\notag\\
&=\int_0^d q_{jj}(t,v) p_{ij}(s,t,u,\dv v)-\frac{\partial p_{ij}}{\partial d}(s,t,u,d),\label{eq:SMGForwardDiagonalLimitResult}
\end{align}
and from (\ref{eq:SMGForwardOffDiagonal}) and (\ref{eq:SMGForwardDiagonalLimitResult}), the theorem follows.
\end{namedproof}

\bibliographystyle{amsplain}

\bibliography{full}

\end{document}